\documentclass[12pt]{amsart}
\usepackage{amsfonts, amssymb,amsmath, amsthm, amsxtra, latexsym,
  mathrsfs, amscd}

\usepackage{amsmath, amssymb, amsthm, times, float}
\usepackage{graphics}
\usepackage{color}
\usepackage{enumitem}   

\newtheorem{thm}{Theorem}[section]

\newtheorem{lemma}[thm]{Lemma}

\newtheorem{rem}[thm]{Remark}

\newcommand{\bc}{{\mathbb C}}

\newcommand{\br}{{\mathbb R}}

\newcommand{\Ga}{\Gamma}

\newcommand{\lo}{\longrightarrow}

\newcommand{\End}{{\rm End}}
\newcommand{\Ker}{{\rm Ker}}

\newcommand{\ad}{{\rm ad}}

\newcommand{\GL}{{\rm \bf GL}}
\newcommand{\tr}{\textnormal{tr}}
\newcommand{\SU}{\mathbf{SU}}

\newcommand{\fk}{\mathfrak{k}}

\newcommand{\fq}{\mathfrak{q}}

\newcommand{\fl}{\mathfrak{l}}
\newcommand{\fm}{\mathfrak{m}}

\newcommand{\fsl}{\mathfrak{sl}}

\newcommand{\fsp}{\mathfrak{sp}}
\newcommand{\fsu}{\mathfrak{su}}
\newcommand{\fso}{\mathfrak{so}}
\newcommand{\fu}{\mathfrak{u}}

\newcommand{\fz}{\mathfrak{z}}

\newcommand{\diag}{\textnormal{diag}}

\newcommand{\fg}{\mathfrak{g}}

\begin{document}
\title[Local rigidity of complex hyperbolic lattice]
{Local rigidity of complex hyperbolic lattices in semisimple Lie groups
}
\author{Inkang Kim and Genkai Zhang}
\date{}
\address{School of Mathematics\\
     KIAS, Hoegiro 85, Dongdaemun-gu\\
     Seoul, 130-722, Korea}
\email{inkang@kias.re.kr}
\address{
Mathematical Sciences, Chalmers University of Technology and
Mathematical Sciences, G\"oteborg University, SE-412 96 G\"oteborg, Sweden}
\email{genkai@chalmers.se}
\begin{abstract}
We show the local rigidity of complex hyperbolic lattices in classical Hermitian
semisimple Lie groups, $SU(p,q), Sp(2n+2,\mathbb R), SO^*(2n+2), SO(2n,2)$. This reproves or generalizes some results in \cite{GM, KKP, Klingler-inv, Pozzetti}.
\end{abstract}

\begin{quote}
\end{quote}

\footnotetext[3]{Research supported by STINT-NRF grant (2011-0031291) and Swedish Research
Council.
I. Kim gratefully acknowledges the partial support
of  grant  (NRF-2014R1A2A2A01005574) and a warm support of IHES  and
Chalmers University of Technology during his stay. G. Zhang would like
to acknowledge and thank the support from KIAS.}
\maketitle

\baselineskip 1.30pc
\section{Introduction}\label{intro}
After the seminal work of A. Weil \cite{Weil,Rag-book}, many pioneering works
about local rigidity of lattices in semisimple Lie groups have been
done by Raghunathan, Matsushima-Murakami, Goldman-Millson and
others.  As in quasi-Fuchsian deformation of Fuchsian groups, one
could expect a possible deformation of a lattice of a semisimple Lie
group $L$ inside a larger Lie group $G\supset L$. Due to
Margulis superrigidity of higher rank semisimple Lie groups, and to
Corlette's superrigidity of $Sp(n,1), F_{4(-20)}$, there is no local
deformation for lattices in such semisimple Lie groups. Hence a
natural interesting problem is to study  lattices of $L= SO(n,1)$ and $SU(n,1)$
in $G \supset L$. A local
rigidity of lattices of $L=SO(3,1)$ and $SO(4,1)$ inside $G=Sp(n,1)$ is
proved in \cite{KP}. A local rigidity of complex hyperbolic uniform lattices $\Gamma$
of $SU(n,1)$ inside $SU(n+1,1)$ is first studied by Goldman-Millson
\cite{GM}. Note that in this case $H^1(\Gamma, \mathfrak{su}(n+1,1))$
decomposes as $H^1(\Gamma,\mathbb R)\oplus H^1(\Gamma, \mathbb
C^{n+1})$. What Goldman-Millson showed is that there is no deformation
coming from $H^1(\Gamma, \mathbb C^{n+1})$. Indeed there is a
deformation coming from $H^1(\Gamma,\mathbb R)$. But this deformation
corresponds to  the less interesting deformation obtained by deforming $\Gamma$ in $U(n,1)$ by a curve of homomorphism into the centralizer $U(1)$ of $SU(n,1)$ in $U(n,1)$.
Henceforth, when we say `locally rigid', we ignore this kind of deformation through the centralizer. A further generalization of this case to quaternionic
K\"ahler manifolds is studied in \cite{KKP, Klingler-inv}.
In those papers
the embeddings of $L=SU(n, 1)$
are in the classical Lie groups and are obtained
by the standard (or natural) embeddings $\rho: L\to
G$
of $L$ in
$G=Sp(n, 1), SU(2n, 2), SO(4n,4)$.

The infinitesimal rigidity result above is partly determined
by the cohomology group $H^1(\Gamma, \fg)$.
It is  known by Raghunathan \cite{Rag-ajm} that the
cohomology group
$H^1(\Gamma, W)$ of $\Gamma$ acting on an irreducible representation
space $W$ of $L=SU(n, 1)$
vanishes unless $W$ is a symmetric tensor
power $S(\mathbb C^{n+1})$ (or its
dual)
 of the standard representation
$\mathbb C^{n+1}$.
In a recent paper  \cite{Klingler-inv}
Klingler proved two results
on local rigidity for uniform lattices of $L$ in $G$
via  $\rho: L\to G$.   To recall his results we denote $M=S(U(n)\times U(1))\subset L$ and $K\subset G$
 the maximal compact subgroups of the respective
semisimple Lie groups, and $\fm, \fl, \fk, \fg$
the corresponding Lie algebras. Let $Z_K(L)$
be the centralizer of $L$ in $K\subset G$ and
 $\fz_{\fk}(\fl)$ its Lie algebra.
The first result  \cite{Klingler-inv} states that
if there exists a non-zero element
$Z\in \fz_{\fk}(\fl)$
 acting on the representation spaces $S^m(V)$
with certain  positivity property (see
(\ref{positivity})
below for a precise statement)
then there is local rigidity. The second result
weakens the assumption by replacing
$Z\in \fz_{\fk}(\fl)$
by the condition
$Z\in \fz_{\fk}(\fm)
$, namely
by $Z$ being in a (generally) larger subspace
$\fz_{\fk}(\fm)\supset
 \fz_{\fk}(\fl)$. The proof of both results uses
a cohomology theory  of polarized real variation of Hodge structures.

The purpose of the present paper is two-fold.
We shall prove  first
 a general local rigidity result by
assuming existence of certain elements $\fz_{\fk}(\fl)$ by using elementary computations for the usual cohomology
of differential forms instead of  polarized real variation of Hodge structures.
Secondly we consider some natural homomorphisms  $\rho: L=SU(n, 1)\to G$
into the classical groups
$SU(p, q), Sp(2n+2, \mathbb R),
SO^*(2n+2), SO(2n, 2)$.
 We explain briefly
how the homomorphisms $\rho:
L\to G$ are constructed.
The first embedding into $SU(p, q)$
is through the diagonal embedding
of $L$ into $SU(nq_0, q_0)$. There are involutions
on the groups $G=
Sp(2n+2, \mathbb R), SO^*(2n+2), SO(2n, 2)$
whose fixed point subgroup is precisely
$U(n, 1)$, namely $U(n, 1)$ is a
symmetric subgroup of $G$ and
$G/U(n, 1)$ is a non-Riemannian symmetric
space.  We examine further
the decomposition of $\fg$ under $L$. In case
that
the symmetric tensor representations
$S^m(\mathbb C^{n+1})$ of $L$ do not appear in the decomposition
the rigidity follows
immediately from the vanishing theorem of
 Raghunathan \cite{Rag-ajm}.
If on the contrary there are such summands
appearing, then we  use  Theorem 1.1
below to prove the local rigidity.

{Throughout this paper, $\Gamma\subset SU(n,1),\ n\geq 2$ is  a uniform lattice.}
We list the  homomorphisms just mentioned:
\begin{enumerate}[label=(\roman*)]
 \item The diagonal homomorphism
of $L=SU(n, 1)$   in $SU(p, q)$, 
\item Satake homomorphism
in $G=Sp(n+1, \mathbb R)$, and 
\item Ihara
homomorphisms into Hermitian Lie groups $SO^*(2n)$ and $SO(2n, 2)$.
\end{enumerate}

Our main results are the following.  The precise notations
will be defined in the next section.
Let $\rho: L=SU(n, 1)\to G$ be a Lie group injective homomorphism.
For simplicity we shall view $L$ as a subgroup of $G$.
Recall that  $Z_K(L)=Z_K(\rho(L))
$
the centralizer of $\rho(L)$
in $K$ and $Z_G(L)$ the centralizer in $G$. Let
$\mathfrak z_{\fk}(\fl)$
and $\mathfrak z_{\fg}(\fl)$ be
the respective Lie algebras. Let $\fg$ be the Lie algebra
of $G$ and $B(\cdot, \cdot)$
the Killing form. The space
$\fg$ is then a real representation of $L$,
  $ \text{Ad}\circ\rho:
 L\mapsto \text{End}(\fg)$,
via $\rho$ and the adjoint representation of $G$ on $\fg$.
As $L$ is real semisimple the space $\fg$
is decomposed
into irreducible representations. Denote $V=
\mathbb C^{n+1}$,  and write its
 decomposition under $U(n)$ as $V=V_1+\mathbb C$,
with $V_1=\mathbb C^n$, the {standard} representation of $U(n)$.
The symmetric powers $S^m(V)$
and $S^m(\mathbb C^{n})$, $m\ge 1$ are irreducible
representations of $L$ and  $U(n)$, respectively.
For simplicity we will denote the representation spaces
  by $S^m$ and  $s^m$.
 We have
 \begin{equation}
   \label{eq:g-1}
\fg=\fg_1 +\fg_0
 \end{equation}
where
 \begin{equation}
   \label{eq:g-2}
\fg_1=\sum_{m} S^m\otimes \mathbb R^{d_m}
 \end{equation}
is the sum of isotypes of $S^m$, and
 $\fg_0$ is of  different isotype from $S^m$.
In particular each space
$ S^m\otimes \mathbb R^{d_m}$
is a representation of $\fl+
\mathfrak z_{\fg}(\fl)$.
Consider
the subspace $s^m\subset S^m$ and  its $M$-isotypes in $
 S^m\otimes \mathbb R^{d_m} $,
 \begin{equation}
   \label{W-m}
W_m:
=s^m\otimes \mathbb R^{d_m} \subset
S^m\otimes \mathbb R^{d_m}.
 \end{equation}
The complex multiplication by $i$ on the real space
$W_m$ will be written simply as usual by $X\to iX$, keeping in mind
that all linear forms involved are real linear forms.  It can also
be written as the Lie algebra action of a center element
$T_0$ of $\fu(n)$; see \S2 below.

\begin{thm}\label{mainthm} { Let $L=SU(n,1),\ n\geq 2$ and $\Gamma\subset L$  a uniform lattice.}
Suppose there exist an $L$-invariant  bilinear form $b$
on $\fg^{\mathbb C}$ and an element $Z\in \mathfrak z_{\fk}(\fl)$
such that
\begin{equation}
  \label{positivity}
-b(\ad(Z)(X),  iX)=-b(Z, [X, iX]) >0, \forall X\in W_{m}, m\ge 1.
\end{equation}
Then $\rho: \Gamma\subset L\to G$ is locally rigid.
\end{thm}

\begin{thm}\label{mainthm2} In the totally geodesic
embeddings (i), (ii), (iii) above,
a {uniform} lattice $\Gamma\subset SU(n,1)\subset G$
is locally rigid.
\end{thm}

The claim for the Satake and Ihara homomorphisms
 answers  partly a question of Pansu posed to us.
Recently Pozzetti \cite[Corollary 1.5]{Pozzetti}  proved also
a local rigidity theorem for the diagonal homomorphism
in (1)  using the geometry of Shilov boundary
of bounded symmetric domains. See also \cite{KM} for a different technique.

We shall present an elementary and independent proof
 of the cases in Theorem \ref{mainthm2}.

We would like to thank Pierre Pansu for some helpful
discussions and for raising several questions. We are grateful
to the referee for some insightful suggestions and comments.

\section{Preliminaries}

We fix first some notation
and convention. If $W=\mathbb C^n$ is a complex vector space we shall denote
$W_{\mathbb R}=\mathbb R^{2n}$ the underlying real vector
space by forgetting the complex structure. 
In particular if $W$ is a complex representation
space of a real reductive Lie group $L$ (acting as
complex linear transformation) the above notation
$W_{\mathbb R}$
makes sense as a real representation
of $L$.

{ A standard representation of a $n\times n$ real or
  complex
matrix group $G$ refers to a standard matrix group
action on the vector space $\mathbb R^n$ or $\mathbb C^n$.}
The complexification
of a real Lie algebra $\fg$
will be denoted by $\fg^{\mathbb C}$.
The Killing form
on $\fg$ and $\fg^{\mathbb C}$
will be denoted by $B(\cdot, \cdot)$.


 Let $V=\mathbb C^{n,1}$
be the space 
$\mathbb C^{n,1}$ equipped with  the
 Hermitian form $(x, y)_J=(Jx, y)$ 
of signature $(n, 1)$, 
where  $J$ is
the  diagonal matrix
$J=\diag(1,\cdots, 1, -1)$ and
$(x, y)$ the standard Hermitian form.
The group $L=SU(n, 1)$ consists of complex linear
transformations on $V$ preserving the form $(x, y)_J$
and of determinant $1$.
The symmetric
space $L/M$, $M=S(U(n)\times U(1))=U(n)\subset L$ being the
maximal compact subgroup, can be realized
as the unit ball $B$ in $\mathbb C^n$,
$B=L/M$ with $z=0$ being the base point.
In  $(n+1)\times (n+1)$-matrix realization of $g\in SU(n, 1)$
the action is given as follows
\begin{equation}
\label{g-action}
g=
\begin{bmatrix}
a&b\\
c&d
\end{bmatrix}
\in L, z\in B \mapsto w=gz=(az+b)(cz+d)^{-1}.
\end{equation}
The Lie algebra
$\fu(n, 1)$
 consists
of matrices $X$ such that
$X^\ast J + JX=0$
and $\fl=\fsu(n, 1)$
the subspaces of elements with trace $0$.
Let $\fl=\fu(n)+\fq$ be the Cartan decomposition
of $\fl$. Here $\fu(n)=s(\fu(n) + \fu(1))$ consists of traceless block-diagonal skew-Hermitian
matrices.
The tangent space $T_{x_0}(B)$ at $x_0$
will be  identified with $\fq=\mathbb C^n$ as a real space.

We fix a central element $T_0$ of the maximal compact
subalgebra $\fu(n)$
$$
T_0=(n+1)^{-1}\sqrt{-1} \text{diag}(1, \cdots, 1, -n).
$$

Let $V_1=\mathbb C^n$ be the
{standard} representation of $U(n)$
and $\det$
be the determinant  representation.
The element $T_0$ defines
the complex structure on $B$
and we have
\begin{equation}
\label{HC-deco}
\fsl(n+1,\bc)=\fsl(n,\bc) + \mathbb C T_0 + \fq^+ + \fq^-
\end{equation}
where $\fq^{\pm}$ are the holomorphic
and anti-holomorphic tangent
spaces at $z=0\in B$. As a representation of $\fu(n)$
we have $\fq^+= V_1\otimes \det$.

The {standard} representation $V$ of $\fl$
under $\fu(n)$ is
$$
V= V_1  \oplus  {\det}^{-1}.
$$
Note that we
have
\begin{equation}
  \label{eq:deco-2}
S^{j}(V)=\oplus_{k=0}^j S^k(V_1)\otimes {\det}^{k-j}, \quad
S^{j}(V')=\oplus_{k=0}^j S^k(V_1')\otimes {\det}^{j-k},
\end{equation}
where  $S^k(V_1)$ is the symmetric tensor of the standard
representation $V_1$ of $\fsu(n)$ and $V'$ is the dual representation of $V$.

%

\section{General rigidity results}

\subsection{Discrete groups and
automorphic cohomologies.
}\label{automorphic}

Let $\Gamma$ be a cocompact lattice of $L=SU(n,1),\ n\geq 2$
and $\rho$ a representation of $\Gamma$ in  a complex vector space
$W$. The infinitesimal deformations of $\rho$ are described by the
cohomology group $H^1(\Gamma, W)
$.
We refer to \cite{Rag-book} for an account
of the underlying theory. 
We have the following result of Raghunathan \cite{Rag-ajm}.

\begin{thm}[Raghunathan]\label{Rag}
Let $\rho: \SU(n,1) \lo \GL(W)$ be a real finite dimensional
irreducible representation of $\SU(n,1), n\geq 2$. Let $\Ga$ be
a cocompact lattice in $SU(n,1)$. Then $H^1(\Ga, W)=0$ except if
$W\simeq  S^j V$ for some $j \geq 0$, where $S^j$ denotes the
$j$-th symmetric power.
\end{thm}

Given a uniform lattice $\Gamma\subset L$
and a holomorphic representation $s$
of the complexification $M^{\mathbb C}=GL(n, \mathbb C)$ 
of $M$ there is also an automorphic
holomorphic bundle over $X=\Gamma\backslash L/M$
and the corresponding
automorphic
cohomology  $H^{p, q}(X, \Gamma, s)$.
In \cite{MM1, MM2} general  relations between
the discrete group cohomology 
associated to a representation of $G$
and the automorphic
cohomology associated to a holomorphic representation of $M^{\mathbb C}$
are studied.
Refined relations between them have
also been further obtained for general
Hermitian Lie groups; see e.g. \cite{Zucker2}.
In  \cite{KKP,KZ, Klingler-inv}  it is proved
for $L=SU(n, 1)$  that the cohomology
 $H^1(\Gamma,  S^j)$ is isomorphic
to an  automorphic $(0, 1)$-cohomology taking values in $s^j$. We
formulate it  in terms of the tangent and
the canonical bundle.
Let $T=T^{(1, 0) }$ be the holomorphic tangent bundle of $X$
and $\mathcal L^{-1}$
be the
line bundle on $X$
defined
 so that $\mathcal L^{-(n+1)}$
is the canonical line bundle $K_X$.
Recall the notation $V=\mathbb C^{n+1}$
and  $V_1=\mathbb C^{n}$, the  representation spaces
of $L$ and $M$.

\begin{thm}
\label{E-S}
Any element $\alpha \in H^1(\Gamma,  S^j(V))$
can be realized as a  $S^j(\mathbb C^n)$-valued
$(0, 1)$-form. This realization induces an isomorphism
$$
H^1(\Gamma,  S^j( V))=
H^{1}(X,  E(S^j(V_1)))
=
H^{(0,1)}(X, S^j TX\otimes \mathcal L^{-j}).
$$
\end{thm}
{For a complete proof, see Theorem 1.4.1 of \cite{Klingler-inv} or Theorem 1.1 of \cite{KZ}.}
\begin{rem}
The cohomology group
$H^1(\Gamma,  S^j(V))$
 can be represented
by harmonic forms with values in a holomorphic
vector bundle over $L/M$ and there is further a decomposition
$H^1(\Gamma,  S^j(V))=
H^{(0, 1)}(\Gamma,  S^j(V
))
+H^{(1, 0)}(\Gamma,  S^j(V))$.
Generally there is an injection of the first part $H^{(0, 1)}(\Gamma,  S^j(V))$
into the automorphic cohomology $H^{1}(X,  E(S^j(V_1)))$; see
\cite{MM1}.
The above theorem implies that the injection is also surjective and
that second term vanishes. In the present paper we shall use
only the injective property of $H^1(\Gamma,  S^j(V))$ into the automorphic
cohomology.
\end{rem}

\subsection{Proof of Theorem 1.1}

\begin{proof} We use similar arguments as
in \cite{KKP} for quaternionic forms combined
with the Eichler-Shimura isomorphism
for the representation $S^m(V)$ of $L$.
Assume that
$\phi\in H^1(\Gamma,
\fg)$ is a non-zero element.
Recall that we have a bilinear form
$[\cdot, \cdot]: H^1(\Gamma,
\fg)
\times
H^1(\Gamma,
\fg)\to
H^2(\Gamma,
\fg)
$,
$[\alpha\otimes X,
\beta\otimes Y,
](u, v)
=(\alpha(u)\beta(v)-
\alpha(v)\beta(u))
[ X,  Y]$. If $\phi$ represents a non-trivial
deformation then we have
$[\phi, \phi]=0$
as an element of $H^2(\Gamma, \fg)$, that
is
$[\phi, \phi]=d\eta$
for some one co-chain $\eta$.
 We shall prove that this leads to a contradiction.
We will freely view $\phi$ and $\eta$ as
differential forms on $G$; see \cite{Rag-book}.
Recall the decomposition
(\ref{eq:g-1}) and the subspaces
$\mathfrak g_1$ and $W_m$ defined in  (\ref{eq:g-2})-
(\ref{W-m}).
It follows first by Theorem 3.1 that $\phi$ takes
values in $\fg_1$.
Theorem 3.2 asserts further that any element in $H^1(\Gamma, S^m(\mathbb
C^{n+1}))$ is a $(0, 1)$-form taking value in the subspace $s^m=S^m(\mathbb
C^n)$. See also \cite{KZ} for an elementary proof. Thus $\phi=\sum_{m}\phi_m$
with $\phi_m$ being $W_m$-valued $(0,1)$-forms.

Now let
 $Z\in \mathfrak z_{\fk}(\fl)$
be  as in the assumption. In particular $Z$ is $L$-invariant.
It defines then a non-zero element
in the  cohomology $H^0(\Gamma, \fg)$,
the subspace of $\Gamma$ invariant elements
in $\fg$ via the action $\text{Ad}\circ\rho$. We consider the bilinear form
$$
\xi(X,Y)=b(Z, [\phi(X),  \phi(Y)])
=b([Z, \phi(X)],   {\phi} (Y) ).
$$

The Lie bracket $[X, Y]$ defines  a $M$-equivariant
map $(\sum_m W_m)
\otimes (\sum_m W_m)
\to \fg$ and we have
$$
[\phi,  {\phi}
 ]
=\sum_{m, m'}
[\phi_m,  {\phi_{m'}}]
$$
and
$$
b(Z, [\phi(X),
{\phi}(Y)
])=\sum_{m, m'} b(Z, [\phi_m(X),  {\phi_{m'}}(Y)])
$$
However for any $m\ne m'$ the spaces $W_m
$
and  $W_{m'}
$ are inequivalent representations
of $M$,  $b$ and $Z$ are $L$-invariant,  thus it follows from Schur's lemma
\cite[6.1]{Hum}
 that
$$
 b(Z, [\phi_m(X), {\phi_{m'}}(Y)])=0
$$
 for any $X, Y$.
Indeed, the bilinear 
form $ b(Z, [\phi_m, \phi_{m'}])$
defines also an $M$-equivariant map from 
$W_{m'}$ to $(W_{m})^\ast=W_{m}$, and it must be zero by Schur's lemma.

 Consequently
$$b(Z, [\phi(X),
{\phi}(Y)
])=\sum_{m} b(Z, [\phi_m(X),  {\phi_{m}}(Y)]).
$$

 Using the fact that the $\phi_m$ are (0,1)-forms and  our assumption,
when we plug in  an orthonormal basis of the underlying real
Euclidean vector space of $\bc^n$, $(E_1 , iE_1
\cdots ,E_n ,iE_n)$, there exists  some $c_x>0$ on each $x\in \Gamma\backslash L/M$ such that
\begin{equation*}
\begin{split}
\frac{\xi\wedge\omega^{n-1}}{\omega^{n}} =\sum_{k=1}^{n}\xi(E_k
,iE_k)&=\sum_{k=1}^n \sum_m b(Z, [\phi_m(E_k),\phi_m(iE_k)])\\
  &=-\sum_{k=1}^n \sum_m b(Z, [\phi_m(E_k),i\phi_m(E_k)])
\geq c_x |\phi|^2_x.
\end{split}
\end{equation*}
Here $\omega$ is a K\"ahler form on $H^n_\bc=L/M$.

{Put $\xi=b(Z, [\phi,\phi])=\lambda\circ [\phi,\phi]$ where $\lambda$ is a linear functional.}
 Therefore
$$\lambda\circ [\phi,\phi]
\wedge \omega^{n-1} \geq c |\phi|^2\omega^n.$$ Here
 $c>0$  since
$\Gamma$ is a uniform lattice.  On the other hand $
[\phi,\phi]=d\eta$ and hence
$$\int_{\Gamma\backslash L/M} \lambda\circ [\phi,\phi]\wedge \omega^{n-1}=0 ,$$ a contradiction.
\end{proof}

\subsection{Proof of Theorems 1.2}
This follows from Theorem 3.1, Theorem \ref{mainthm} and  Lemmas 4.1-4.3.

\section{Embeddings of the unit ball in classical symmetric spaces}\label{examples}

In this section we construct
natural homomorphisms
of
$L=SU(n, 1)$
in classical Lie groups $G=SU(p, q),
Sp(p,  \mathbb R), SO^*(2p), SO(2, p)$, and we prove
the relevant decomposition of the Lie
algebra $\fg$ under $\fsu(n, 1)$. We examine
the condition in the statement of Theorem 1.1. The Killing form
for these classical Lie algebra will be fixed once for all
as $B(X, Y)=\tr XY$.
In the $ SO^*(2p)$ and $ SO(2,p)$ cases, we will show that $S^m(V)$ factor does not appear
in the Lie algebra decomposition, hence by Theorem
\ref{Rag}, the local rigidity follows. 
(Certain decompositions
of the complex Lie algebras here can also be done as in \cite{Eastwood-Wolf}.)

\subsection{
The diagonal homomorphism
of $L=SU(n, 1)$
in  $G=SU(p, q)$
}

We consider first the case 
  $G=SU(nq_0, q_0)$. Let $
W:=V\otimes \mathbb C^{q_0}
=\mathbb C^{n, 1}\otimes \mathbb C^{q_0}
$ be
equipped with the Hermitian inner product
$(\cdot, \cdot)_W:
=(\cdot, \cdot)_J\otimes (\cdot, \cdot)
$, where $ (\cdot, \cdot)$ is the standard Hermitian inner product
on $\mathbb C^{q_0}$. 
 The  groups
$L=SU(n, 1)$ and $U(n, 1)$  are diagonally embedded in $SU(W)=SU(nq_0, q_0)$
and respectively in $U(W)=U(nq_0, q_0)$
via $g\to g\otimes I_{\mathbb C^{q_0}}$,
and we 
fix this realization. To find  $L$-invariant forms as in Theorem 1.1
it is natural to consider the centralizer of $L$ in $SU(nq_0, q_0)
$. Indeed
the group $U(q_0)$ is also diagonally embedded in $G$
via $h\to I_{\mathbb C^{n, 1}}\otimes h$, and  commutes with
 $U(n, 1)$ in $SU(nq_0, q_0)$. (The pair
$(U(n, 1), U(q_0))$ is an example of the  so-called Howe
dual pairs \cite{Howe}.) 
{The space of $L$-invariant bilinear forms on $\mathfrak{su}(nq_0, q_0)$
forms then a representation space of $U(q_0)$ and it
is therefore conceptually clear to
consider the  decomposition of 
$\mathfrak{su}(nq_0, q_0)$ first  under $L \times U(q_0)$.}

Consider the Lie algebra 
$\fu(nq_0, q_0)$ of $U(W)=U(nq_0, q_0)$.
Let $\mathcal H=\mathcal H(q_0)=
\{Y\in M_{q_0, q_0}; Y^\ast=Y \}$ be the space of $q_0\times q_0$-Hermitian
matrices $Y$ viewed as a representation space of $U(q_0)$, namely
$h\in U(q_0): Y\mapsto hYh^\ast$; $\mathcal H$ is identified
with the Lie algebra $\fu(q_o)=i\mathcal H$ of $U(q_0)$ as a representation
space.
 It is immediate that
the  Lie algebra 
$\fu(nq_0, q_0)$ under the standard action of $U(n, 1)\times U(q_0)$ is
decomposed
as 
$$
\fu(nq_0, q_0)
=\mathfrak{u}(n, 1) \otimes \mathcal H.
$$
Indeed any  element
of the form $X\otimes Y\in \mathfrak{u}(n, 1)
\otimes \mathcal H$ is clearly an element in
$\fu(nq_0, q_0)$, since
\begin{equation*}
  \begin{split}
(X\otimes Y (v_1\otimes u_1),
v_2\otimes u_2)_W&=(Xv_1, v_2)_J (Yu_1, u_2)=
-(v_1, Xv_2)_J (u_1, Yu_2)\\
&=
-(v_1\otimes u_1,
X\otimes Y (v_2\otimes u_2))_W,
  \end{split}
\end{equation*}
and the whole space is generated by elements
of this form  by counting the dimensions.

We have now  $\mathfrak u(n, 1)=\fl + i\mathbb R I_{C^{n, 1}}$, 
$\mathcal H=\mathbb RI_{\mathbb C^{q_0}} 
+\mathcal H_0$
where $\mathcal H_0$ is the subspace of trace free elements.
By 
taking the trace free part  and observing that $\tr(X\otimes Y)=
(\tr X)(\tr Y)$ 
we find 
$$\fsu(nq_0, q_0)=\fl \otimes \mathcal H +
i I_{C^{n, 1}}
\otimes \mathcal H_0.
$$
To see how the Lie algebra $\fl$ is realized 
as a subalgebra
we write the above formula as
\begin{equation}
\label{tr-fre}
\begin{split}
\fsu(nq_0, q_0) &=(\mathfrak l \otimes I_{\mathbb C^{q_0}} +
iI_{\mathbb C^{n, 1}}\otimes \mathcal H_0
) +
\mathfrak{l}\otimes \mathcal H_0 
\\
&=(\mathfrak l + \mathfrak{su}
(q_0))
+\mathfrak{l}\otimes \mathcal H_0 
\end{split}
\end{equation}
under the action of  $L\times U(q_0)$. Here we have used 
the above identification of  $L$ and $U(q_0)$ 
as well as their Lie algebras, so that $\fl=
\mathfrak l \otimes I_{\mathbb C^{q_0}}$,
$\mathfrak{su}(q_0)=iI_{\mathbb C^{n, 1}}\otimes \mathcal H_0$.

We consider further  the  homomorphism of $SU(n, 1)$
in $G=SU(p, q)$ via the above embedding $SU(n, 1)\to SU(nq_0, q_0)$
and the natural inclusion  $SU(nq_0, q_0)\subset
SU(p, q)
$ for $p\ge nq_0, q\ge q_0$ and $p+q >(n+1)q_0$.
More precisely let $p_1=p-nq_0, q_1=q-q_0$ and let
$\mathbb C^{p, q}
=\mathbb C^{nq_0, q_0} \oplus
\mathbb C^{p_1, q_1}$
 be equipped with the indefinite Hermitian
form defined by the matrix
$$
\text{diag}(J_1, J_2), \, \text{where}\, 
J_1=\text{diag}(I_{nq_0}, -I_{q_0}),
\, J_2=\text{diag}(I_{p_1}, -I_{q_1}).
$$
The group
$SU(p, q)$ is then $SU(\mathbb C^{p, q})$.
The Lie algebra elements  
in $\mathfrak{g}$ will
be written as $2\times 2$-block matrices
under the above decomposition of 
$\mathbb C^{p, q}$, and they are of the form
$$
\begin{bmatrix} a&b\\
-J_2b^*J_1
&d \end{bmatrix},
$$
with $a\in 
\mathfrak u(nq_0, q_0),
d\in \mathfrak u(p_1, q_1),
$ $\tr a+ \tr d=0$.
Let
$$
E:=\begin{bmatrix} i(p_1+q_1)I_{nq_0 +q_0}
 &0\\
0 &
-i (n+1)q_0I_{p_1+q_1}
\end{bmatrix}\in \fg,
$$
be the sum of two central elements  in $\mathfrak u(nq_0, q_o)$
and  respectively in  $\mathfrak{u}(p_1, q_1)$,
which is further in $\mathfrak{z}_{\fk}(\fl)$ according to our notation in Section \ref{intro}.
Thus we have
\begin{equation}
  \label{step-1}
\mathfrak{g}=
\mathfrak{su}(nq_0, q_0)
\oplus
\mathfrak{su}(p_1, q_1) 
\oplus \mathbb RE \oplus
\left(\mathbb C^{nq_0, q_0}\otimes_{\mathbb C}
(\mathbb C^{p_1, q_1})'\right)_{\mathbb R}
\end{equation}
as a representation space of 
$\mathfrak {su}(nq_0, q_0)
+ \mathfrak{su}(p_1, q_1)
$ with 
$\mathbb C^{nq_0, q_0}$ the standard action
of  $\mathfrak {su}(nq_0, q_0)$,
 $(\mathbb C^{p_1, q_1})'$
the dual action of 
$\mathfrak{su}(p_1, q_1)$ and
$\mathbb R E$
the trivial representation of  
$\mathfrak {su}(nq_0, q_0)
+ \mathfrak{su}(p_1, q_1)$.

{The last factor $\mathbb C^{nq_0, q_0}\otimes_{\mathbb C}
(\mathbb C^{p_1, q_1})'$ can be seen as follows. The matrices $A\in U(nq_0,q_0)$ and $ B\in U(p_1,q_1)$ act
on $$\begin{bmatrix} 0&b\\
-J_2b^*J_1
&0 \end{bmatrix},
$$ by $$\begin{bmatrix} A &0\\
0
&B \end{bmatrix} \begin{bmatrix} 0&b\\
-J_2b^*J_1
&0 \end{bmatrix}\begin{bmatrix} A^{-1} &0\\
0
&B^{-1} \end{bmatrix}=\begin{bmatrix} 0 & AbB^{-1}\\
*
&0 \end{bmatrix}.
$$}

Now the subgroup 
 $L\times SU(q_0)\subset
SU(nq_0, q_0)\subset G$
acts on 
$\mathbb C^{nq_0, q_0}$ as $\mathbb C^{n, 1}\otimes \bc^{q_0}$
and thus on  the last
summand
in   (\ref{step-1})  as
$$
\mathbb C^{nq_0, q_0}\otimes
(\mathbb C^{p_1, q_1})'
=\mathbb C^{n, 1}\otimes \mathbb C^{q_0}
\otimes (\mathbb C^{p_1, q_1})'.
$$
The first summand is treated in   (\ref{tr-fre}), and 
the formula   (\ref{step-1}) 
now reads 
\begin{equation}
  \label{eq:triple}
  \begin{split}
\mathfrak{g}&=
\left(\mathfrak{l} + \mathfrak{su}(q_0)\right)
\oplus (\mathfrak{l}\otimes \mathcal H_0)
\oplus \mathfrak{su}(p_1, q_1) \oplus \mathbb RE 
 \\
&\qquad  \oplus (\mathbb C^{n, 1}\otimes_{\mathbb C}\mathbb
  C^{q_0}\otimes_{\mathbb C}(\mathbb C^{p_1, q_1})')_{\mathbb R}.
      \end{split}
\end{equation}
 {Note here that the diagonal embedding of $SU(n,1)$ in $G$ is via $L\times SU(q_0)\subset SU(nq_0,q_0)\subset G$ and hence it commutes with
$SU(p_1,q_1)$. Therefore the above decomposition is under the diagonal embedding of $SU(n,1)$. }
Therefore the standard representation $\bc^{n,1}$ of $L$
appears with multiplicity $q_0(p_1+q_1)$.
Using the notation in Section \ref{intro} we have
$$
W_1=\mathbb C^{n}\otimes \mathbb C^{q_0}
\otimes (\mathbb C^{p_1, q_1})',
$$
consisting of matrices 
$$
X(b)=\begin{bmatrix} 0&b\\
-J_2 b^*J_1 &0 \end{bmatrix},
$$ with 
$$
b=\begin{bmatrix} b_{11}&b_{12}\\
0&0\end{bmatrix},
$$
$b_{11}$ of size $nq_0\times p_1$
and  $b_{12}$ of size 
$nq_0\times q_1$.

Assume first $p_1=p-nq_0 >0$. We let
$$
Z=\text{diag}(i \alpha I_{(n+1)q_0}, i\beta I_{p_1}, i\gamma I_{q_1})
\in \mathfrak z_{\mathfrak k}(\mathfrak l),
$$
where $\alpha, \beta, \gamma$ 
are real numbers subjected to $(n+1)q_0\alpha+p_1\beta+q_1\gamma=0$
so that $\tr Z=0$,  and will be chosen afterwards. 
In the computation below we shall
suppress the index $m$ of  $I_{m}$.
 
We compute the action of $Z$ on $X$
and find
$$
[Z, X(b)]
=\begin{bmatrix} 0&v\\
-J_2v^*J_1 &0 \end{bmatrix},
$$
with 
$$
v =i \alpha I
\begin{bmatrix} b_{11}&b_{12}\\
0&0\end{bmatrix}
- i\begin{bmatrix} b_{11}&b_{12}\\
0&0\end{bmatrix}  \text{diag}(\beta I, \gamma I)
=\begin{bmatrix} i(\alpha -\beta) b_{11}& i (\alpha -\gamma) b_{12}\\
0 &0 \end{bmatrix}.
$$

We check now the condition in Theorem \ref{mainthm}.
Replacing $b$ by $ib$ we have
$$
X(ib)=\begin{bmatrix} 0&ib\\
iJ_2 b^*J_1 &0 \end{bmatrix}.
$$ 
The bilinear form $B([Z, X(b)], X(ib))$
in question is
$$
B([Z, X(b)], X(ib))
=\tr [Z, X(b)]X(ib)
=-2(\alpha -\beta)
 \tr b_{11} b_{11}^*  
-2(\gamma-\alpha)\tr b_{12} b_{12}^*
$$
by a direct computation. We choose
now $\beta=-1$, $\gamma$ any real number such that
$$
\frac{p_1}{(n+1)q_0 +q_1} <\gamma  <
\frac{(n+1)q_0 +p_1}{q_1}
$$
(which clearly exists) and 
$$
\alpha=\frac{p_1 -\gamma q_1}{(n+1)q_0}
$$
Then we have  indeed $
(n+1)q_0\alpha+p_1\beta+q_1\gamma=0$,
and
$$
\gamma>\alpha >\beta
$$
so that 
$B([Z, X(b)], X(ib))$ is negative definite.

If $p_1=p-nq_0 =0$ then $q_1=q-q_0>0$. We
replace $Z$ above by the following matrix
$$
Z=\text{diag}(-i\frac{q_1}{(n+1)q_0}I_{(n+1)q_0}, iI_{q_1})
$$
The same computation as above shows that we still have
$$
B([Z, X], iX)=-  2\tr  b^* b <0
$$
whenever $b\ne 0$.

Summarizing we have
\begin{lemma}\label{po-su}  Consider
the diagonal embedding
of $L=SU(n, 1)$
in $G=SU(p, q)$ via $L\times U(q_0)\subset G$ for $p\ge nq_0, q\ge q_0$
\begin{enumerate}
\item Suppose $p=nq_0, q=q_0$.
The  Lie algebra $\fg$
is decomposed under $\mathfrak{su}(n, 1) + \mathfrak{su}(q)$
as
$$
\fg=\fsu(nq_0, q_0)=
\left(
\mathfrak{su}(n, 1)
\otimes \mathcal H({q})
\right)
\oplus
\left(
I_{\mathbb C^{n, 1}}
\otimes \fsu(q)\right)=(\mathfrak l + \mathfrak{su}
(q))
\oplus (\mathfrak{l}\otimes \mathcal H_0(q))
$$
where $\mathcal H(q)
$ is the space of $q\times q$-Hermitian
matrices viewed as a representation space of $\mathfrak{su}(q)$ and $\mathcal H_0(q)$ is the trace free part.
Thus no symmetric tensor $S^mV$ of $
\mathfrak{su}(n, 1)
$ appears in
the decomposition {under the diagonal embedding of $SU(n,1)$ in $G$.}
\item Suppose $p+q>nq_0+ q_0$.  We have
$\mathfrak{g}$ is decomposed under
$\mathfrak l+\mathfrak {u}(q_0)
+\mathfrak{su}(p_1, q_1)$ as
in  (\ref{eq:triple}).
There exists an element $Z\in \mathfrak z_{\fk}(\fl)$
such that the positivity condition in Theorem 1.1 holds.
\end{enumerate}
\end{lemma}

\subsection{
The homomorphism of $\fl=\fsu(n, 1)$
into  $\fg=\fso(2n,  2),
\fso^*(2n+2)$
}

These homomorphisms have been studied
by Ihara \cite{Ihara} in terms of root systems.
We shall need more precise formulation.

We shall realize $U(n, 1)$ as symmetric
subgroup of
$G=SO^*(2n+2)$
and $SO(2n, 2)$. The relevant
decompositions of the Lie algebras $\fg=
\fso^*(2n+2)$ and $\fso(2n, 2)$ under
$\fu(n, 1)$ are parallell
to the Cartan decomposition of the noncompact
symmetric space
$SO^*(2n+2)/U(n+1)$
and the compact symmetric space  $SO(2n+2)/U(n+1)$. For completeness we provide details here.

We consider first  the group $SO(2n, 2)$.
Let $\mathbb R^{2n+2}=\mathbb C^{n+1}$ as a real space
be equipped with the real indefinite form $\Re(x, y)_J=\Re (Jy)^* x$,
where $(x, y)_J$ is the Hermitian form on $\mathbb C^{n+1}$
of signature $(n, 1)$ in \S2.1. The Lie group $G=SO(2n, 2)$
is the connected component of the identity of the group $O(2n, 2)$
defined by the real indefinite form. The symmetric
space of $G=SO(2n, 2)$ is the Lie ball,
a complex structure being determined by fixing
an element of the Lie algebra of $\fso(2)$.
In particular $SU(n, 1)$ is a subgroup of $SO(2n, 2)$, and the corresponding
embedding of the unit ball into the Lie ball
is holomorphic by choosing a consistent complex structure on the ball.
The multiplication
by $i$, $z\to iz$ on $\mathbb R^{2n+2}$
 induces an involution  $\theta$ on the Lie algebra
$ \mathfrak{so}(2n, 2),
\theta^2=id$. Thus
$ \mathfrak{so}(2n, 2)= \ker(\theta-1)\oplus
\ker(\theta+1)
$ and it is clear that
$\ker(\theta-1)=\mathfrak{u}(n, 1)$, the subspace
of  $\mathbb C$-linear transformations $T: z\mapsto Tz$
in $ \mathfrak{so}(2n, 2)$.
The subspace $\ker(\theta+1)$ consists of  anti-$\mathbb C$-linear
transformations of the form $\alpha_A:z\mapsto A\bar z$  which are anti-symmetric with
respect
to the real form $\Re(z, w)_J$, i.e., $$\Re (A\bar z, w)_J= -\Re (z, A \bar w)_J.$$
That
is $\Re (Jw)^* A \bar z=-\Re (JA\bar w)^*z
=-\Re (JA\bar w)^t \bar z$, namely
$JA=-A^tJ$. Putting $\tilde A=JA$ we have
$\tilde A^t =-\tilde A$ and the action of an element
$T\in \mathfrak{u}(n, 1)$ on $\alpha_A\in \text{ker}(\theta+1)$ is
$$[T, \alpha_A]:z\mapsto T\alpha_A z- \alpha_AT z=TA\bar z- A(\overline{Tz})=(TA-A\bar T)\bar z.$$ Hence
$A\mapsto TA-A\bar T
$. In terms of $\tilde A$
this is $\tilde A \mapsto J(TA-A\bar T)$, which
is
 $$
J TA-JA\bar T
=-T^*JA -JA\bar T =-(T^* \tilde A + \tilde A \bar T)
$$
using  $T\in \mathfrak{u}(n, 1)$, $JT=-T^* J$. Namely
$\mathfrak{u}(n, 1)$ acts on $ \text{ker}(\theta+1)$
via the (real) representation on the space $V\wedge V$ of anti-symmetric $(n+1)\times (n+1)$ complex matrices.

We consider now the group
 $G=SO^*(2n+2)$.  Let us  fix its realization first.
Let
$$
j=\begin{bmatrix} 0 & J\\
-J& 0\end{bmatrix},
$$
where $J=J_{n, 1}=\diag(I_n, -1)$ is the matrix defining $U(n, 1)$ in \S2.1.
The group $G$
can be  realized  as
$$
G=\{g\in GL(2n+2, \mathbb C);  g^*jg=j, g^tg=I_{2n+2}, \det g=1\}.
$$
Note that $G$ here is a different realization
from the one in \cite{He1}
but
it will be convenient for our purpose.
The group
$SO^*(2n+2)$ is realized in
\cite[Ch. X, \S2]{He1}  as
$$
G_1=
\{g\in GL(2n+2, \mathbb C);  g^*j_1g=j_1, g^tg=I_{2n+2}, \det g=1\}.
$$
where
$$
j_1=\begin{bmatrix} 0 & I_{n+1}\\
-I_{n+1}& 0\end{bmatrix}.
$$
It is clear that $j_1$ and  $j$ are conjugate, $c^{-1}jc=j_1$, by an
element $c$  in $O(2n+2)$, which interchanges $e_{n+1}$ and $e_{2n+2}$, so that the two groups $G$
and $G_1$ are isomorphic by the map $g\in G_1 \mapsto cgc^{-1}\in G$.

Let $\tau: g\to jgj^{-1}$ be the involution on $
GL(2n+2, \mathbb C)$. Then $\tau$
maps $G$ to $G$ and is an involution on $G$
and  $\fg$.
We claim that the set of fixed point
is the group $U(n, 1)$
and respectively $\mathfrak{u}(n,1)$.   To be more precise
let $\fg=\fg_{+}
+\fg_{-}
=\Ker (\tau-1)\oplus
\Ker (\tau+1)$
 be the decomposition
of $\fg$ under $\tau$. Then
$\fg_{+} $ is  a symmetric subalgebra of $\fg$ and by elementary matrix
computations we find that the Lie algebra
$\fg_{+} $ consists of real skew symmetric matrices
$$
X=\begin{bmatrix} A & B\\
C& D\end{bmatrix}, \quad \bar X=X, \quad X^t=-X,
$$
satisfying
$$
AJ=JD, \quad BJ =JB^t .
$$
It is precisely the Lie algebra $\fu(n, 1)$ under the
identification of $X$ with the complex matrix $A-iBJ$,
i.e., $\fg_+=\fu(n, 1)$.
The subspace $\fg_{-} $ consists of complex matrices of the form
$iX$ where $X$ are real skew symmetric matrices
$$\quad X=\begin{bmatrix} A & B\\
C& D\end{bmatrix}, \quad \bar X=X, \quad X^t=-X,
$$
with
$$
AJ=-JD, \quad BJ=-JB^t .
$$
We identify $\fg_-$ with the space $
V\wedge V=\mathbb C^{n, 1}\wedge \mathbb C^{n, 1}
$
of skew-symmetric matrices
via the map $iX\mapsto A+iBJ$. The Lie algebra action of $\fg_+
=\fu(n, 1)
$ on $\fg_-$ is also the tensor product
action of $\fu(n, 1)$ on $V\wedge V$, by an elementary
matrix computation of the  Lie bracket in $\fg$.
 Summarizing we obtain the following

\begin{lemma}\label{po-so} Under
the above realizations of
${SU}(n, 1)$ as a subgroup in $G=
SO(2n, 2)$ and $SO^*(2n+2)$ the  Lie algebra $\fg$
 has decomposition
$\fg=\mathfrak{u}(n, 1)
\oplus(V\wedge V)_{\mathbb R}$  under $\fl=\mathfrak{su}(n, 1)$.
In particular the symmetric tensors $S^mV$ do not
appear in the decomposition.
\end{lemma}

\subsection{
The Lie algebra  $\fg=\fsp(n+1,  \mathbb R)$
and  subalgebra
$\fl=\fsu(n, 1)$
}

View $V=\mathbb C^{n+1}$ as a real space $V=\mathbb R^{2n+2}$
with the complex structure $v\to iv$ and the
Hermitian form $(x, y)_J$ as in \S2.
Consider the following
symplectic form
$$
\omega (x, y)=\Im(J x, y)
$$
on $V$.
We let $Sp(n+1,\mathbb R)
=:Sp(n+1,\omega)$ be the symplectic
group defined by $\omega$.  Clearly $U(n, 1)$
is a subgroup of $Sp(n+1,\mathbb R)$ and
the unit ball  $U(n, 1)/U(n)\times U(1)$
is embedded holomorphically into the Siegel domain
$Sp(n+1,\mathbb R)/U(n+1)$.
Symplectic embeddings of Hermitian
symmetric spaces have been studied systematically
by Satake \cite{Satake}.

The complex multiplication $i$ on $V$ defines
an involution $A\mapsto i^{-1}Ai$ on all real  linear
transformations $A\in \End(\mathbb R^{2n+2})$,
giving a decomposition
of $\End(\mathbb R^{2n+2})=M_{n+1, n+1}(\mathbb C)
\oplus {M_{n+1, n+1}(\mathbb C)}^c$
with $M_{n+1, n+1}(\mathbb C)$ consisting  of complex linear transformations
$T: x\to Tx$
of $\mathbb C^{n+1}=\mathbb R^{2n+2}$,
and ${M_{n+1, n+1}(\mathbb C)}^c$
of the conjugate complex linear
transformations $T^c: x\to T\bar x$
where $T\in M_{n+1, n+1}(\mathbb C)$.
Restricting the involution to
the Lie algebra $sp(n+1,\mathbb R)$
we find
$$
 \mathfrak{sp}(n+1,\mathbb R) =
\mathfrak{u}(n, 1)
\oplus \mathfrak r
$$
where
$$
\mathfrak r
=\{T^c; (JT)^t =JT\},
$$
where $A^t$ is the transpose of $A$. Indeed $\mathfrak r$ consists of
complex conjugate transformations
$T^c: x\to T\bar x$ in $
 \mathfrak{sp}(n+1,\mathbb R)$,
i.e.,
$
\omega (T^c x, y) + \omega(x, T^cy)
=0$, $x, y\in V$; equivalently
\begin{eqnarray*}
0&=&\Im (JT\bar x, y) + \Im(x, JT\bar y)
=\Im y^\ast JT \bar x + \Im (JT \bar y)^\ast x\\
&=&\Im y^\ast JT \bar x + \Im y^t \overline{(JT )^t} x
=\Im y^\ast JT \bar x - \Im y^\ast {(JT )^t} \bar x,
\end{eqnarray*}
implying $(JT)^t=JT$. As a representation space of $\mathfrak{u}(n, 1) $,
$\mathfrak r$ is identified
with the symmetric tensor power
$S^2(V)
$
of the standard representation $V$,
  concretely via
$T^c \to JT\in S^2(V)$.
We have further, via this  identification,
$$ \mathfrak{sp}(n+1,\mathbb R)
=\mathfrak{su}(n, 1) \oplus\mathbb R Z_0
\oplus \mathfrak r=
\mathfrak{su}(n, 1)
\oplus\mathbb R Z_0
\oplus (S^2(V))_{\mathbb R}
$$
where the element $Z_0$ generates the centralizer
of $\mathfrak{su}(n, 1)$  in $\mathfrak{u}(n, 1)$
defined as the multiplication by $i$ on $V_\mathbb R$.
 The Lie
bracket on $\mathfrak r$ is given,  by the definition,
$$
[X^c, Y^c]= U,  \quad U=X\bar Y - Y\bar X\in \mathfrak{u}(n, 1).
$$
In particular
we have
$$
[X^c, iX^c]
= -2i X\bar X.
$$
The Killing form
on $\fsp(n+ 1,\br)$ is $B(X, Y)=\tr_{\mathbb R} XY$
by our normalization.
Now if $X, Y\in S^2(V)$
is in the last component
$S^2(V_1)\subset
S^2(V)$ in
  (\ref{eq:deco-2}),
then $(JX)^t=X^t=JX=X$ since $J=I$ on $V_1$.
Hence
we have
$$
B(Z_0, [X^c, i X^c])=\tr_{\mathbb R} i ([X^c, i X^c])=\tr_{\mathbb R} i(-2i X\bar X)=2\tr_{\mathbb R}X\bar X=2\tr_{\mathbb R}X X^*
=4\tr_{\mathbb C} XX^\ast,
$$
and is positive definite.
Summarizing we have
\begin{lemma}\label{po} The Lie algebra
$ \mathfrak{sp}(n+1,\mathbb R)$ has a decomposition
$\fg= \mathfrak{sp}(n+1,\mathbb R) =
\mathfrak{su}(n, 1)
\oplus \mathbb R Z_0
\oplus (S^2(V))_{\mathbb R}=
\fl
\oplus
\mathbb R Z_0
\oplus (S^2(V))_{\mathbb R}
$ under $\fl$.
The $\fl$-invariant
bilinear skew-symmetric form $(X^c, Y^c)\to
B(Z_0, [X^c, Y^c])$
on $S^2(V)$ induces
a positive definite quadratic form
$X^c\rightarrow
B(Z_0, [X^c, iX^c])$
 on the subspace $
S^2(V_1)\subset S^2(V)$.
\end{lemma}

\begin{rem}
Consider  a natural homomorphism
of $L=SU(2, 1)$ into the
exceptional Lie group $G=F_{4(-20)}$
inducing a totally geodesic embedding
of the complex unit ball $L/M$ in $\mathbb C^2$
into  $G/K$ realized as the unit ball in the octonian
$\mathbb O^2$ (see e.g. \cite[Theorem 2.2]{Johnson-2}).
We can show in this case that
$
\mathfrak z_{\fk}(\fm)=\mathfrak z_{\fk}(\fl)=\mathfrak{su}(3)$. But we observe that for any element
$Z\in \mathfrak z_{\fk}(\fm)
=\mathfrak z_{\fk}(\fl)$ the  positivity condition
(\ref{positivity})
is never fulfilled. 
It is thus an open question
to know whether the local rigidity holds in this case.
\end{rem}

\end{document}